\documentclass{article}
\usepackage{amsfonts}
\usepackage{amsmath}

\input{tcilatex}
\begin{document}

\author{I. Bairamov$^{1}$\thinspace \thinspace \thinspace and K. Bayramoglu$%
^{2}$ \\
$^{1}${\small Department of Mathematics Izmir University of Economics,
Izmir, Turkey. }\\
{\small E-mail: ismihan.bayramoglu@ieu.edu.tr}\\
$^{2}${\small Department of Statistics, Middle East Technical University,
Ankara,Turkey.}\\
{\small E-mail: konul@metu.edu.tr}}
\date{\thinspace \thinspace \thinspace\ \ \ \ \ \ }
\title{Baker- Lin-Huang type Bivariate distributions based on order
statistics }
\maketitle

\begin{abstract}
Baker (2008) introduced a new class of bivariate distributions based on
distributions of order statistics from two independent samples of size $n$.
Lin-Huang (2010) discovered an important property of Baker's distribution
and showed that the Pearson's correlation coefficient for this distribution
converges to maximum attainable value, i.e. the correlation coefficient of
the Frech\'{e}t upper bound, as $n$ increases to infinity. Bairamov and
Bayramoglu (2011) investigated a new class of bivariate distributions
constructed by using Baker's model and distributions of order statistics
from dependent \ random variables, allowing high correlation than that of
Baker's distribution. \ In this paper \ a new class of Baker's type
bivariate distributions with high correlation are constructed on the base of
distributions of order statistics by using an arbitrary continuous copula
instead of the product copula.

\thinspace \thinspace \thinspace


\thinspace \thinspace \thinspace
\end{abstract}

\section{Introduction}

Huang and Kotz (1999) introduced new modifications of classical
Farlie-Gumbel-Morgenstern (FGM) distribution introducing additional
parameters. The new Huang-Kotz FGM distributions allow correlation higher
than the classical FGM and because of the simple analytical form aroused
interest of many researchers. Last years there appeared many papers dealing
with the modifications of FGM distribution allowing high correlation. For a
related works on this subject see \ Bairamov et al. (2001), Amblard and
Girard (2002), Bairamov and Kotz (2002), (2003), Fisher and Klein (2007)
among others. \ Baker (1998) use a novel approach connected with the FGM
distribution and introduce a new class of bivariate distributions based on
the distributions of order statistics. Baker (2008) considers independent
random variables $X_{i},$ $Y_{i}$, from two univariate distributions with
distribution functions (c.d.f.) $\ F_{X}$ and $F_{Y},$ respectively. The
corresponding probability density functions (p.d.f.) are $f_{X}$ and $f_{Y}. 
$ Let $U_{1}=\min (X_{1},X_{2})$ and $U_{2}=\max (X_{1},X_{2}),$ $V_{1}=\min
(Y_{1},Y_{2})$ and $V_{2}=\max (Y_{1},Y_{2}).$ To obtain positive
correlation Baker randomly chooses either the pair $U_{1},V_{1}$ or $%
U_{2},V_{2}$. The random numbers are now positively correlated, but the
marginal distributions are unchanged, because the random choice of either of
two order statistics from a distribution gives a random variable from that
distribution. To obtain a negative correlation, either $U_{1}$,$V_{2}$ or $%
U_{2},V_{1}$ chosen. The bivariate distribution of a randomly chosen pair of
order statistics is 
\begin{equation*}
1/2\times \{F_{X}^{2:2}(x)F_{Y}^{2:2}(y)+F_{X}^{1:2}(x)F_{Y}^{1.2}(y)\}
\end{equation*}%
and on choosing either a pair \ of order statistics with probability $q$ or
the original independent random variables $X,Y$ with probability $1-q,$ the
bivariate distribution function is 
\begin{eqnarray*}
H(x,y)
&=&(1-q)F_{X}(x)F_{Y}(y)+(q/2)%
\{F_{X}^{2:2}(x)F_{Y}^{2:2}(y)+F_{X}^{1:2}(x)F_{Y}^{1.2}(y)\} \\
&=&F_{X}(x)F_{Y}(y)\{1+q(1-F_{X}(x))(1-F_{Y}(y)\}.
\end{eqnarray*}%
In general let $X_{1},X_{2},...,X_{n}$ and $Y_{1},Y_{2},...Y_{n}$ are
independent and identically distributed (i.i.d.) random variables with
distribution functions (d.f.) $F_{X}$ and $F_{Y},$ respectively. Let $%
X_{k:n} $ and $Y_{k:n},$ $k=1,2,...,n$ be corresponding order statistics and 
$F_{X}^{k:n}(x)=P\{X_{k:n}\leq x\},$ $F_{Y}^{k:n}(y)=P\{Y_{k:n}\leq y\}.$
Baker's bivariate distribution function is now defined as 
\begin{eqnarray}
H_{+}^{(n)}(x,y) &=&\frac{1}{n}\sum_{k=1}^{n}F_{X}^{k:n}(x)F_{Y}^{k:n}(y)
\label{b1} \\
H_{-}^{(n)}(x,y) &=&\frac{1}{n}%
\sum_{k=1}^{n}F_{X}^{k:n}(x)F_{Y}^{n-k+1:n}(y).  \label{b2}
\end{eqnarray}%
For Baker's bivariate distribution $H_{+}^{(n)}(x,y)$ with exponential
marginals $F_{X}(x)=F_{Y}(x)=1-e^{-x},x>0,$ the Pearson's correlation
coefficient is $\rho _{n}=1-\frac{1}{n}\sum_{k=1}^{n}\frac{1}{k},$ which
increases monotonely to 1.

As a generalization of (\ref{b1}) and (\ref{b2}) Baker also introduces%
\begin{equation}
H_{r}^{(n)}(x,y)=\dsum\limits_{k=1}^{n}%
\sum_{l=1}^{n}r_{kl}F_{X}^{k:n}(x)F_{Y}^{l:n}(y),  \label{b3}
\end{equation}%
where $r_{kl}\geq 0$ and $\sum_{k=1}^{n}r_{kl}=\sum_{l=1}^{n}r_{kl}=\frac{1}{%
n},$ for all $k,l=1,2,...,n.$ Lin and Huang (2010) proved that (\ref{b3})
does not contain members with a correlation higher than that of (\ref{b1}).
That is why the best bivariate distribution with higher positive correlation
among the members (\ref{b3}) is (\ref{b1}). \ Similar consideration holds
true for the negative correlation. Lin and Huang (2010) investigate the
conditions under which the correlation for (\ref{b1}) converges to the
limit. In particular, they show that if either (i) $X\geq b,Y\geq c$ a.s.
for some $b,c\in 
\mathbb{R}
$ and $E(X_{k:n})\geq F_{X}^{-1}(\frac{k-1}{n})$ and $E(Y_{k:n})\geq
F_{Y}^{-1}(\frac{k-1}{n})$ for all ($k,n),$ or (ii) $X\leq b,Y\leq c$ a.s.
for some $b,c\in 
\mathbb{R}
$ and $E(X_{k:n})\leq F_{X}^{-1}(\frac{k}{n})$ and $E(Y_{k:n})\leq
F_{Y}^{-1}(\frac{k}{n})$ for all ($k,n)$ then $\underset{n\rightarrow \infty 
}{\lim }\rho _{n}=\rho ^{\ast },$ where $\rho ^{\ast }$ is the correlation
coefficient of the Fr\'{e}chet-Hoeffding upper bound $H_{+}(x,y)=\min
(F_{X}(x),F_{Y}(y))$ (see Frech\'{e}et (1940)). The results presented in the
paper of Lin and Huang (2010) makes Baker's distribution attractive.

Recently, Bairamov and Bayramoglu (2011) observed that if in the Baker's
model instead of independent random variables one uses the dependent random
variables $(X,Y)$ with positive quadrant dependent (PQD) joint distribution
function $F(x,y),$ then the correlation increases and for negative quadrant
dependent $F(x,y),$ decreases. More precisely, let ($%
X_{1},Y_{1}),(X_{2},Y_{2}),...,(X_{n},Y_{n})$ be a bivariate sample with
joint distribution function $F(x,y)=C(F_{X}(x),F_{Y}(y))$. \ Bairamov and
Bayramoglu (2011) consider the following bivariate distribution functions
constructed on the basis of the Baker's idea: 
\begin{eqnarray*}
K_{+}^{(n)}(x,y) &=&\frac{1}{n}\sum_{r=1}^{n}P\{X_{r:n}\leq x,Y_{r:n}\leq y\}
\\
K_{-}^{(n)}(x,y) &=&\frac{1}{n}\sum_{r=1}^{n}P\{X_{r:n}\leq
x,Y_{n-r+1:n}\leq y\},
\end{eqnarray*}%
where $X_{i:n}$ and $Y_{j:n}$ are the $i$th and $j$th order statistics
constructed on the basis of bivariate observations ($X_{i},Y_{i}),$ ($%
i=1,2,...,n)$ with joint distribution function $F(x,y)=P\{X_{i}\leq
x,Y_{i}\leq y\}$ and marginal distribution functions $F_{X}(x)=F(x,\infty ),$
$F_{Y}(y)=F(\infty ,y)$ so that $X_{1:n}\leq X_{2:n}\leq \cdots \leq
X_{n:n}, $ $Y_{1:n}\leq Y_{2:n}\leq \cdots \leq Y_{n:n}.$ \ The joint d.f.
of $X_{r:n} $ and $Y_{s:n}$ is given in David (1981) (see also Arnold et al.
(1992)) as%
\begin{eqnarray}
P\{X_{r:n} &\leq &x,Y_{s:n}\leq y\}  \notag \\
&=&\dsum\limits_{i=r}^{n}\sum_{j=s}^{n}%
\sum_{k=a}^{b}c(n,k;i,j)p_{11}^{k}p_{12}^{i-k}p_{21}^{j-k}p_{22}^{n-i-j+k},
\label{x1}
\end{eqnarray}%
where 
\begin{equation}
c(n,k;i,j)=\frac{n!}{k!(i-k)!(j-k)!(n-i-j+k)!},\text{ }a=\max (0,i+j-n),%
\text{ }b=\min (i,j)  \label{xx1}
\end{equation}%
and 
\begin{eqnarray}
p_{11} &=&F(x,y)  \notag \\
p_{12} &=&F_{X}(x)-F(x,y)  \notag \\
p_{21} &=&F_{Y}(y)-F(x,y)  \notag \\
p_{22} &=&1-F_{X}(x)-F_{Y}(y)+F(x,y).  \label{xxx1}
\end{eqnarray}%
Then%
\begin{eqnarray}
K_{+}^{(n)}(x,y) &=&\frac{1}{n}\dsum\limits_{r=1}^{n}\dsum\limits_{i=r}^{n}%
\sum_{j=r}^{n}%
\sum_{k=a}^{b}c(n,k;i,j)p_{11}^{k}p_{12}^{i-k}p_{21}^{j-k}p_{22}^{n-i-j+k}
\label{x2} \\
K_{-}^{(n)}(x,y) &=&\frac{1}{n}\dsum\limits_{r=1}^{n}\dsum\limits_{i=r}^{n}%
\sum_{j=n-r+1}^{n}%
\sum_{k=a}^{b}c(n,k;i,j)p_{11}^{k}p_{12}^{i-k}p_{21}^{j-k}p_{22}^{n-i-j+k}.
\label{x3}
\end{eqnarray}

It is clear that the marginal distributions of $\ K_{+}^{(n)}(x,y)$ and $%
K_{-}^{(n)}(x,y)$ are again $F_{X}(x)$ and $F_{Y}(y),$ i.e. $%
K_{+}^{(n)}(x,\infty )=F_{X}(x)$ and $K_{-}^{(n)}(\infty ,y)=F_{Y}(y).$ It
is shown that for a PQD joint distribution function $F(x,y)$ the positive
correlation of $K_{+}^{(n)}(x,y)$ is higher than that of $H_{+}^{(n)}(x,y)$
and for NQD $F(x,y)$ and the negative correlation of $K_{-}^{(n)}(x,y)$ is
smaller than that of $H_{-}^{(n)}(x,y).$

In this note we consider a new class of bivariate distribution functions
using Baker's \ construction and considering \ any copula $C(u,v)$ instead
of product copula $C(u,v)=\Pi (u,v)=uv.$ It follows that for this new class
of distributions if $C(u,v)$ is PQD, i.e. $C(u,v)\geq uv,$ for all $(u,v)\in
\lbrack 0,1]^{2},$ then the Pearson's correlation coefficient is higher than
that of Baker's distribution. Similarly, if the copula is NQD, i.e. $%
C(u,v)\leq uv,$ for all $(u,v)\in \lbrack 0,1]^{2}$ then then the Pearson's
correlation coefficient is smaller than that of Baker's distribution. All
modifications constructed on the Baker's idea we  call as Baker-Lin-Huang
Type distributions. \ The distribution (\ref{x2}), (\ref{x3}) \ we call
Baker's Type I BB and the new distributions introduced in this paper Baker's
Type II BB distributions.

\section{New bivariate Baker's Type II BB distributions based on a copula
approach}

\ \ Let $X_{1},X_{2},...,X_{n}$ and $Y_{1},Y_{2},...Y_{n}$ \ be i.i.d.
random variables with d.f.'s $F_{X}$ and $F_{Y},$ respectively. Let $X_{k:n}$
and $Y_{k:n},$ $k=1,2,...,n$ be corresponding order statistics and $%
F_{X}^{k:n}(x)=P\{X_{k:n}\leq x\},$ $F_{Y}^{k:n}(y)=P\{Y_{k:n}\leq y\}.$

Recall that a two-dimensional copula is a function $C(x,y)$ from $%
[0,1]^{2}=[0,1]\times \lbrack 0,1]$ to $[0,1]$ with the properties:

1. $C(x,0)=0=C(0,y),$ $C(x,1)=x$ and $C(1,y)=y;$

2. For every $x_{1},x_{2},y_{1},y_{2}$ such that $0\leq x_{1}<x_{2}\leq 1$
and $0\leq y_{1}<y_{2}\leq 1$ 
\begin{equation*}
C(x_{2},y_{2})-C(x_{2},y_{1})-C(x_{1},y_{2})+C(x_{1},y_{1})\geq 0.
\end{equation*}%
According to Sklar's Theorem if $F_{X,Y}(x,y)$ is a joint distribution
function with continuous marginal distributions $F_{X}(x)$ and $F_{Y}(y),$
then there exists a unique copula $C$ such that $%
F(x,y)=C(F_{X}(x),F_{Y}(y)). $ Theory and applications of copulas are well
documented in Nelsen (2005). Let $C(u,v)$ be any copula. Consider 
\begin{eqnarray}
G_{+}^{(n)}(x,y) &=&\frac{1}{n}\sum_{k=1}^{n}C(F_{X}^{k:n}(x),F_{Y}^{k:n}(y))
\label{m1} \\
G_{-}^{(n)}(x,y) &=&\frac{1}{n}%
\sum_{k=1}^{n}C(F_{X}^{k:n}(x),F_{Y}^{n-k+1:n}(y)).  \label{m2}
\end{eqnarray}%
From the properties of a copula it follows that the marginal distributions
of $G_{+}^{(n)}(x,y)$ and $G_{-}^{(n)}(x,y)$ are $F_{X}$ and $F_{Y},$
respectively. In fact, 
\begin{eqnarray*}
G_{+}^{(n)}(x,\infty ) &=&\frac{1}{n}\sum_{k=1}^{n}C(F_{X}^{k:n}(x),1)=\frac{%
1}{n}\sum_{k=1}^{n}F_{X}^{k:n}(x)=F_{X}(x) \\
G_{+}^{(n)}(\infty ,y) &=&\frac{1}{n}\sum_{k=1}^{n}C(1,F_{Y}^{k:n}(y))=\frac{%
1}{n}\sum_{k=1}^{n}F_{Y}^{k:n}(y)=F_{Y}(y).
\end{eqnarray*}%
Similarly, $G_{-}^{(n)}(x,\infty )=F_{X}(x)$ and $G_{-}^{(n)}(\infty
,y)=F_{Y}(y).$ Hereafter, we will denote by $\rho _{H}$ $\ $the Pearson's
correlation coefficient between any random variables $X$ and $\ Y$ with
joint distribution function $H(x,y).$ It is clear that if $C(u,v)=\Pi
(u,v)=uv$ then $G_{+}^{(n)}(x,y)=H_{+}^{(n)}(x,y)$ and $%
G_{-}^{(n)}(x,y)=H_{-}^{(n)}(x,y).$ Since $C(F_{X}^{k:n}(x),F_{Y}^{k:n}(y))$
is a bivariate c.d.f. (with marginals $F_{X}^{k:n}(x)$ and $F_{Y}^{k:n}(y))$
then $G_{+}^{(n)}(x,y)$ is obviously a bivariate c.d.f. \ as a convex
combination of bivariate c.d.f.'s. The copula used in construction (\ref{m1}%
) and (\ref{m2}) will be called the "kernel" copula for Baker's Type II BB
distribution.

\ \ 

\textbf{Theorem 1. }If $C(u,v)$ is PQD then $\rho _{G_{+}^{(n)}}\geq \rho
_{H_{+}^{(n)}}$ and if $C(u,v)$ is NQD then $\rho _{G_{-}^{(n)}}\leq \rho
_{H_{-}^{(n)}}.$

\textbf{Proof.} \ Since $C(u,v)\geq uv$ for all $(u,v)\in \lbrack 0,1]^{2},$
then $C(F_{X}^{k:n}(x),F_{Y}^{k:n}(y))\geq F_{X}^{k:n}(x)F_{Y}^{k:n}(y)$ for
all $(x,y)\in 
\mathbb{R}
^{2}.$ From the Hoeffding's formula (see Hoeffding (1940)) for correlation
coefficient one has 
\begin{eqnarray}
\rho _{G_{+}^{(n)}} &=&\frac{1}{\sqrt{Var(x)Var(Y)}}\dint\limits_{-\infty
}^{\infty }\dint\limits_{-\infty }^{\infty
}[G_{+}^{(n)}(x,y)-F_{X}(x)F_{Y}(y)]dxdy  \notag \\
&=&\frac{1}{\sqrt{Var(x)Var(Y)}}\dint\limits_{-\infty }^{\infty
}\dint\limits_{-\infty }^{\infty }[\frac{1}{n}%
\sum_{k=1}^{n}C(F_{X}^{k:n}(x),F_{Y}^{k:n}(y))  \notag \\
&&-F_{X}(x)F_{Y}(y)]dxdy  \label{m3} \\
\rho _{H_{+}^{(n)}} &=&\frac{1}{\sqrt{Var(x)Var(Y)}}\dint\limits_{-\infty
}^{\infty }\dint\limits_{-\infty }^{\infty
}[H_{+}^{(n)}(x,y)-F_{X}(x)F_{Y}(y)]dxdy  \notag \\
&=&\frac{1}{\sqrt{Var(x)Var(Y)}}\dint\limits_{-\infty }^{\infty
}\dint\limits_{-\infty }^{\infty }[\frac{1}{n}%
\sum_{k=1}^{n}F_{X}^{k:n}(x)F_{Y}^{k:n}(y)  \notag \\
&&-F_{X}(x)F_{Y}(y)]dxdy  \label{m4}
\end{eqnarray}%
and $\rho _{G_{+}^{(n)}}\geq \rho _{H_{+}^{(n)}}.$ Similarly, $\rho
_{G_{-}^{(n)}}\leq \rho _{H_{-}^{(n)}}.$

\bigskip\ \ 

\textbf{Example 1. Baker's Type II BB distributions with FGM \ "kernel"
copula.} \ Let 
\begin{equation}
C(u,v)=uv(1+\alpha (1-u)(1-v)),(u,v)\in \lbrack 0,1]^{2},-1\leq \alpha \leq
1.  \label{FGM}
\end{equation}%
Consider 
\begin{eqnarray}
G_{+}^{(n)}(x,y) &=&\frac{1}{n}\sum_{k=1}^{n}F_{X}^{k:n}(x)F_{Y}^{k:n}(y)(1+
\label{m5} \\
\alpha (1-F_{X}^{k:n}(x))(1-F_{Y}^{k:n}(y)),\text{ }0 &\leq &\alpha \leq 1. 
\notag \\
G_{-}^{(n)}(x,y) &=&\frac{1}{n}%
\sum_{k=1}^{n}F_{X}^{k:n}(x)F_{Y}^{n-k+1:n}(y)(1+  \notag \\
\alpha (1-F_{X}^{k:n}(x))(1-F_{Y}^{n-k+1:n}(y)),\text{ }-1 &\leq &\alpha
\leq 0.  \label{m6}
\end{eqnarray}

Let $F_{X}(x)=x,0\leq x\leq 1$ and $F_{Y}(y)=y,0\leq y\leq 1.$ Since the FGM
copula (\ref{FGM}) is PQD for $\alpha \geq 0$ and is NQD for $\alpha \leq 0,$
\ \ then $\rho _{G_{+}^{(n)}}\geq \rho _{H_{+}^{(n)}}$ for $\alpha \geq 0,$ $%
\ \rho _{G_{-}^{(n)}}\leq \rho _{H_{-}^{(n)}},$ for $\alpha \leq 0.$ \ It is
clear that if $\alpha =0$ then $G_{+}^{(n)}(x,y)=H_{+}^{(n)}(x,y)$ and $%
G_{-}^{(n)}(x,y)=H_{-}^{(n)}(x,y).$ In the following table we present some
numerical values of \ $\rho _{G_{+}^{(n)}},\rho _{K_{+}^{(n)}},\rho
_{H_{+}^{(n)}},$ $\rho _{G_{-}^{(n)}},\rho _{K_{-}^{(n)}}$ and $\rho
_{H_{-}^{(n)}}$for $Uniform(0,1)$ marginals. The numerical calculations are
made in MATLAB which is one of the commonly accepted packages for coding
mathematical models since it has built-in functions for probability
distributions and allows probabilistic and mathematical operations.

\ \ \bigskip

\begin{equation*}
\begin{tabular}{|l|l|l|l|l|l|l|l|l|}
\hline
$n$ & $2$ & $4$ & $6$ & $8$ & $10$ & $12$ & $15$ & $20$ \\ \hline
$\rho _{G_{+}^{(n)}}$ & $0.5467$ & $0.7258$ & $0.8039$ & $0.8475$ & $0.8753$
& $0.8945$ & $0.9144$ & $0.9348$ \\ \hline
$\rho _{K_{+}^{(n)}}$ & $0.5133$ & $0.6915$ & $0.7761$ & $0.8247$ & $0.8561$
& $0.8779$ & $0.9006$ & $0.9241$ \\ \hline
$\rho _{H_{+}^{(n)}}$ & $0.3333$ & $0.6000$ & $0.7143$ & $0.7778$ & $0.8182$
& $0.8462$ & $0.8750$ & $0.9048$ \\ \hline
$\rho _{G_{-}^{(n)}}$ & $-0.5467$ & $-0.7258$ & $-0.8039$ & $-0.8475$ & $%
-0.8753$ & $-0.8945$ & $-0.9144$ & $-0.9348$ \\ \hline
$\rho _{K_{-}^{(n)}}$ & $-0.5133$ & $-0.6915$ & $-0.7761$ & $-0.8247$ & $%
-0.8561$ & $-0.8779$ & $-0.9006$ & $-0.9241$ \\ \hline
$\rho _{H_{-}^{(n)}}$ & $-0.3333$ & $-0.6000$ & $-0.7143$ & $-0.7778$ & $%
-0.8182$ & $-0.8462$ & $-0.8750$ & $-0.9048$ \\ \hline
\end{tabular}%
\end{equation*}%
\begin{eqnarray*}
&&Table\text{ }1.\text{ Correlation coefficients }\rho _{G_{+}^{(n)}}\text{
and }\rho _{G_{-}^{(n)}}\text{ with } \\
F(x,y) &=&xy(1+\alpha (1-x)(1-y)),\alpha =1,F_{X}(x)=x,F_{Y}(y)=y,0\leq
x,y\leq 1\text{ } \\
&&\text{and }\rho _{K_{+}^{(n)}}\text{, }\rho _{K_{-}^{(n)}}\text{ and }\rho
_{H_{+}^{(n)}}\text{, }\rho _{H_{-}^{(n)}}\text{ with }Uniform(0,1)\text{
marginals}.
\end{eqnarray*}

\bigskip

\bigskip\ \ \ \ \ 

It can be observed from the Table 1 that $\ \rho _{G_{-}^{(n)}}\leq \rho
_{K_{-}^{(n)}}\leq \rho _{H_{-}^{(n)}}\leq \rho _{H_{+}^{(n)}}\leq \rho
_{K_{+}^{(n)}}\leq \rho _{G_{+}^{(n)}}.$

\bigskip In the following table the values of the correlation coefficients
for $\rho _{G_{+}^{(n)}},\rho _{K_{+}^{(n)}}\rho _{H_{+}^{(n)}},$ $\ \rho
_{G_{-}^{(n)}},\rho _{K_{-}^{(n)}}$ and $\rho _{H_{-}^{(n)}}$ for $%
Uniform(0,1)$ and $Exponential(1)$ marginal distributions.

\begin{equation*}
\begin{tabular}{|l|l|l|l|l|l|l|l|l|}
\hline
$n$ & $2$ & $4$ & $6$ & $8$ & $10$ & $12$ & $15$ & $20$ \\ \hline
$\rho _{G_{+}^{(n)}}$ & $0.4811$ & $0.6343$ & $0.7001$ & $0.7367$ & $0.7600$
& $0.7762$ & $0.7929$ & $0.8102$ \\ \hline
$\rho _{K_{+}^{(n)}}$ & $0.4426$ & $0.5951$ & $0.6682$ & $0.7105$ & $0.7379$
& $0.7572$ & $0.7772$ & $0.7980$ \\ \hline
$\rho _{H_{+}^{(n)}}$ & $0.2886$ & $0.5196$ & $0.6185$ & $0.6735$ & $0.7085$
& $0.7327$ & $0.7577$ & $0.7835$ \\ \hline
$\rho _{G_{-}^{(n)}}$ & $-0.4811$ & $-0.6343$ & $-0.7001$ & $-0.7367$ & $%
-0.7600$ & $-0.7762$ & $-0.7929$ & $-0.8102$ \\ \hline
$\rho _{K_{-}^{(n)}}$ & $-0.4426$ & $-0.5951$ & $-0.6682$ & $-0.7105$ & $%
-0.7379$ & $-0.7572$ & $-0.7772$ & $-0.7980$ \\ \hline
$\rho _{H_{-}^{(n)}}$ & $-0.2886$ & $-0.5196$ & $-0.6185$ & $-0.6735$ & $%
-0.7085$ & $-0.7327$ & $-0.7577$ & $-0.7835$ \\ \hline
\end{tabular}%
\end{equation*}%
\begin{eqnarray*}
&&Table\text{ }2.\text{ Correlation coefficients }\rho _{G_{+}^{(n)}}\text{
and }\rho _{G_{-}^{(n)}}\text{ with } \\
F(x,y) &=&xy(1+\alpha (1-x)(1-y)),\alpha =1,F_{X}(x)=x,0\leq x\leq
1,F_{Y}(y)=1-\exp (-y),y\geq 0\text{ } \\
&&\text{and }\rho _{K_{+}^{(n)}}\text{, }\rho _{K_{-}^{(n)}}\text{ and }\rho
_{H_{+}^{(n)}}\text{, }\rho _{H_{-}^{(n)}}\text{ with the same marginals. }
\end{eqnarray*}

\bigskip\ \ \ \ \ \ \ 

Again, \bigskip from the Table 2 we have $\rho _{G_{-}^{(n)}}\leq \rho
_{K_{-}^{(n)}}\leq \rho _{H_{-}^{(n)}}\leq \rho _{H_{+}^{(n)}}\leq \rho
_{K_{+}^{(n)}}\leq \rho _{G_{+}^{(n)}}.$

\section{Copula representation of Baker's Type II BB distribution}

Denote by $\Pi (t,s)=ts,$ $(t,s)\in \lbrack 0,1]^{2},$ a product copula. Let
($X_{i},Y_{i}),i=1,2,..n$ be a bivariate sample with joint distribution
function $F(x,y)=C(F_{X}(x),F_{Y}(y))$. \ Consider the joint distribution of
order statistics $X_{r:n}$ and $Y_{s:n}$ given in (\ref{x1}):%
\begin{eqnarray}
F_{X_{r:n},Y_{s:n}}(x,y) &=&P\{X_{r:n}\leq x,Y_{s:n}\leq y\}  \notag \\
&=&\dsum\limits_{i=r}^{n}\sum_{j=s}^{n}%
\sum_{k=a}^{b}c(n,k;i,j)p_{11}^{k}p_{12}^{i-k}p_{21}^{j-k}p_{22}^{n-i-j+k},
\label{xx2}
\end{eqnarray}%
where $c(n,k;i,j),$ $p_{11},p_{12},p_{21},p_{22}$ are given in (\ref{xx1}),(%
\ref{xxx1}). \ The copula of vivariate distribution $F_{X_{r:n},Y_{s:n}}(x,y)
$ presents an interest. \ Denote this copula as $C_{r,s:n}(t,s),$ then 
\begin{equation}
F_{X_{r:n},Y_{r:n}}(x,y)=C_{r,s:n}(F_{X_{r:n}}(x),F_{Y_{s:n}}(y)).
\label{x4}
\end{equation}

It is well known (David (1981)) 
\begin{eqnarray*}
F_{X_{r:n}}(x) &=&\dsum\limits_{i=r}^{n}\binom{n}{i}%
F_{X}^{i}(x)(1-F_{X}(x))^{n-i} \\
&=&\frac{1}{B(r,n-r+1)}\dint%
\limits_{0}^{F_{X}(x)}u^{r-1}(1-u)^{n-r}du=I_{r,n-r+1}(F_{X}(x))
\end{eqnarray*}%
and 
\begin{eqnarray*}
F_{Y_{s:n}}(y) &=&\dsum\limits_{i=s}^{n}\binom{n}{i}%
F_{Y}^{i}(y)(1-F_{Y}(y))^{n-i} \\
&=&\frac{1}{B(s,n-s+1)}\dint%
\limits_{0}^{F_{Y}(y)}u^{s-1}(1-u)^{n-s}du=I_{s,n-s+1}(F_{Y}(y)),
\end{eqnarray*}%
where $I_{a,b}(p)=\frac{1}{B(a,b)}\dint\limits_{0}^{p}u^{a-1}(1-u)^{b-1}du$
is an incomplete Beta function. Denote by $I_{a,b}^{-1}(p)$ the inverse of $%
I_{a,b}(p).$ Let \ $F_{X_{r:n}}(x)=I_{r,n-r+1}(F_{X}(x))=t$ and $%
F_{Y_{s:n}}(y)=I_{s,n-s+1}(F_{Y}(y))=s.$ Then 
\begin{equation}
x=F_{X}^{-1}(I_{r,n-r+1}^{-1}(F_{X}(t)))\text{ \ and \ \ }%
y=F_{Y}^{-1}(I_{s,n-s+1}^{-1}(F_{Y}(s))).  \label{x6}
\end{equation}%
From (\ref{x4}) and (\ref{x6}) one has 
\begin{equation}
C_{r,s:n}(t,s)=F_{X_{r:n},Y_{r:n}}\left(
F_{X}^{-1}(I_{r,n-r+1}^{-1}(t)),F_{Y}^{-1}(I_{s,n-s+1}^{-1}(s))\right) .
\label{x7}
\end{equation}%
Therefore a copula of joint distribution of order statistics $X_{r:n}$ and $%
Y_{s:n}$ is $C_{r,s:n}(t,s)$ given in (\ref{x7}). In a special case if $%
r=s=n $ one has 
\begin{eqnarray}
C_{n,n:n}(t,s) &=&F^{n}\left(
F_{X}^{-1}(I_{n,1}^{-1}(t)),F_{Y}^{-1}(I_{n,1}^{-1}(s))\right)  \label{x8} \\
&=&F^{n}(F_{X}^{-1}(t^{1/n}),F_{Y}^{-1}(s^{1/n})),  \notag
\end{eqnarray}%
since $I_{n,1}(t)=t^{n}$ and $I_{n,1}^{-1}(v)=v^{1/n}.$

Analogously, one obtains the copula of joint distribution of order
statistics $X_{1:n}$ and $Y_{1:n}$ as 
\begin{eqnarray}
C_{1,n:n}(t,s) &=&t+s-1  \notag \\
&&+[(1-t)^{1/n}+(1-s)^{1/n}-1  \label{x9} \\
&&+F(F_{X}^{-1}(1-(1-t)^{1/n}),F_{Y}(1-(1-s)^{1/n}))]^{n}  \notag
\end{eqnarray}%
by noting that the joint distribution function of $X_{1:n}$ and $Y_{1:n}$ is 
\begin{eqnarray}
F_{X_{1:n},Y_{1:n}}(x,y) &=&P\{X_{1:n}\leq x,Y_{1:n}\leq y\}  \notag \\
&=&(1-(1-F_{X}(x))^{n})+(1-(1-F_{Y}(y))^{n})-1+\bar{F}^{n}(x,y)  \notag \\
&=&(1-(1-F_{X}(x))^{n})+(1-(1-F_{Y}(y))^{n})-1  \notag \\
&&+(1-F_{X}(x)-F_{Y}(y)+F(x,y))^{n}  \notag \\
&=&C_{1,n:n}\left( 1-(1-F_{X}(x)))^{n},1-(1-F_{Y}(y))^{n}\right) .
\label{x10}
\end{eqnarray}

\subsection{The Baker's Type II BB distributions with "kernel" copula of
bivariate FGM extreme order statistics\ }

If the underlying distribution is classical FGM, i.e. 
\begin{equation}
F(x,y)=F_{X}(x)F_{Y}(y)(1+\alpha (1-F_{X}(x))(1-F_{Y}(y))),-1\leq \alpha
\leq 1.  \label{xx10}
\end{equation}%
one obtains from (\ref{x8})%
\begin{equation}
C_{n:n}(t,s)=ts\left[ 1+\alpha (1-t^{1/n})(1-s^{1/n})\right] ^{n}
\label{mm7}
\end{equation}%
and from (\ref{x9}) 
\begin{eqnarray}
C_{1:n}(t,s) &=&t+s-1+\left[ (1-t)^{1/n}+(1-s)^{1/n}-1\right.  \notag \\
&&+(1-(1-t)^{1/n})(1-(1-s)^{1/n})  \notag \\
&&\left. \times \{1+\alpha (1-t)^{1/n}(1-s)^{1/n}\}\right] ^{n}.  \label{mm8}
\end{eqnarray}

\bigskip It follows that $\underset{n\rightarrow \infty }{\lim }$ $%
C_{n:n}(t,s)=ts=\Pi (t,s)$ $\ $and $\underset{n\rightarrow \infty }{\lim }%
C_{1:n}(t,s)=\Pi (t,s),$ since $\underset{n\rightarrow \infty }{\underset{%
n\rightarrow \infty }{\lim }\left[ 1+\alpha (1-t^{1/n})(1-s^{1/n}) 
_{\substack{  \\ }}\right] ^{n}=1\text{ and }}$%
\begin{eqnarray*}
&&\underset{n\rightarrow \infty }{\lim }\left[ (1-t)^{1/n}+(1-s)^{1/n}-1%
\right. \\
&&\left. +(1-(1-t)^{1/n})(1-(1-s)^{1/n})\{1+\alpha (1-t)^{1/n}(1-s)^{1/n}\} 
\right] ^{n} \\
&=&(1-t)(1-s),
\end{eqnarray*}%
where $\Pi (t,s)$ is a product copula of independent random variables. This
means that if the joint distribution of $(X,Y)$ is FGM given in (\ref{xx10}%
), then the extreme order statistics $X_{n:n}$ and $Y_{n:n}$ are
asymptotically independent. So are the order statistics $X_{1:n}$ and $%
Y_{1:n}.$ Despite of this fact, the Baker's type BB distribution constructed
on the base of the copula (\ref{mm7}) and (\ref{mm8}) has enough large
correlation. Indeed, consider Baker's type BB distribution (\ref{m1}) with
the "kernel" copula (\ref{mm7}) and with the "kernel" copula (\ref{mm8}):%
\begin{eqnarray*}
\check{G}_{+}^{(n)}(x,y) &=&\frac{1}{n}%
\sum_{k=1}^{n}C_{n,n:n}(F_{X}^{k:n}(x),F_{Y}^{k:n}(y)) \\
&=&\frac{1}{n}\sum_{k=1}^{n}F_{X}^{k:n}(x)F_{Y}^{k:n}(y)\left[ 1+\alpha (1-%
\left[ F_{X}^{k:n}(x)\right] ^{1/n})(1-\left[ F_{Y}^{k:n}(y)\right] ^{1/n}) 
_{\substack{  \\ }}\right] ^{n} \\
\hat{G}_{+}^{(n)}(x,y) &=&\frac{1}{n}%
\sum_{k=1}^{n}C_{1,n:n}(F_{X}^{k:n}(x),F_{Y}^{k:n}(y)) \\
&=&\frac{1}{n}\sum_{k=1}^{n}\left\{ F_{X}^{k:n}(x)+F_{Y}^{k:n}(y)-1+\left[
(1-F_{X}^{k:n}(x))^{1/n}+(1-F_{Y}^{k:n}(y))^{1/n}-1_{\substack{  \\ }}%
\right. \right. \\
&&+(1-(1-F_{X}^{k:n}(x))^{1/n})(1-(1-F_{Y}^{k:n}(y))^{1/n}) \\
&&\left. \left. \times \{1+\alpha
(1-F_{X}^{k:n}(x))^{1/n}(1-F_{Y}^{k:n}(y))^{1/n}\}\right] ^{n}\right\} .
\end{eqnarray*}%
For example, the correlation coefficient of $\check{G}_{+}^{(n)}(x,y)$ for $%
Uniform(0,1)$ marginals is $\rho _{\hat{G}_{+}^{(n)}}=0.4985$ for $n=2.$%
\bigskip

\subsection{The Baker's type BB distributions with "kernel" copula of
bivariate Gumbel's extreme order statistics}

Let 
\begin{equation}
C(u,v)=\frac{uv}{1+u-uv}.  \label{h1}
\end{equation}%
(\ref{h1}) is the copula of the Gumbel's bivariate logistic distribution
function 
\begin{equation*}
H_{X,Y}(x,y)=(1+e^{-x}+e^{-y})^{-1}\text{, }x\geq 0,y\geq 0\text{ }
\end{equation*}%
with standard logistic marginal distributions $H_{X}(x)=(1+e^{-x})^{-1}$ and 
$H_{Y}(y)=(1+e^{-y})^{-1},$ $\ x\geq 0,y\geq 0.$ (Gumbel (1961), see also
Nelsen (2005), page 28). Using%
\begin{equation*}
F(x,y)=\frac{F_{X}(x)F_{Y}(y)}{F_{X}(x)+F_{Y}(y)-F_{X}(x)F_{Y}(y)},
\end{equation*}%
one obtains from (\ref{x8}) 
\begin{equation}
C_{n:n}(t,s)=\frac{ts}{(t^{1/n}+s^{1/n}-t^{1/n}s^{1/n})^{n}}  \label{kk7}
\end{equation}%
and from (\ref{x9}) 
\begin{eqnarray}
C_{1:n}(t,s) &=&t+s-1+  \notag \\
&&+\left[ _{\substack{  \\ }}(1-t)^{1/n}+(1-s)^{1/n}-1\right.  \label{kk8} \\
&&\left. +\frac{(1-(1-t)^{1/n})(1-(1-s)^{1/n})}{%
2-(1-t)^{1/n}-(1-s)^{1/n}-(1-(1-t)^{1/n})(1-(1-s)^{1/n})}\right. ^{n}. 
\notag
\end{eqnarray}

\bigskip It is seen that 
\begin{equation*}
\lim_{n\rightarrow \infty }C_{n:n}(t,s)=\lim_{n\rightarrow \infty }\frac{ts}{%
(t^{1/n}+s^{1/n}-t^{1/n}s^{1/n})^{n}}=ts=\Pi (t,s).
\end{equation*}

\section{Joint distribution of bivariate order statistics for Frech\'{e}t
upper bound copula}

Let ($X_{1},Y_{1}),(X_{2},Y_{2}),...,(X_{n},Y_{n})$ be a bivariate sample
with joint distribution function $F(x,y)=C(F_{X}(x),F_{Y}(y))$. \ In this
section our aim is firstly, to investigate the joint distribution function
of bivariate order statistics $(X_{r:n},Y_{r:n})$ for a copula with maximal
correlation, i.e. the Frech\'{e}et upper bound. We are interested then in
distribution function 
\begin{equation*}
K_{+}^{(n)}(x,y)=\frac{1}{n}\dsum\limits_{r=1}^{n}\dsum\limits_{i=r}^{n}%
\sum_{j=r}^{n}%
\sum_{k=a}^{b}c(n,k;i,j)p_{11}^{k}p_{12}^{i-k}p_{21}^{j-k}p_{22}^{n-i-j+k},
\end{equation*}%
in the case \ where the marginal distributions are uniform and $C(u,v)=\min
(u,v).$ Recall that the coefficients $c(n,k;i,j)$ and $%
p_{11},p_{12},p_{21},p_{22}$ are given in (\ref{xxx1}),(\ref{xx1}).
Secondly, we consider a distribution introduced in (\ref{m1}) with uniform
marginals, i.e. 
\begin{eqnarray*}
G_{+}^{(n)}(x,y) &=&\frac{1}{n}\sum_{k=1}^{n}C(F_{X}^{k:n}(x),F_{Y}^{k:n}(y))
\\
&=&\frac{1}{n}\sum_{r=1}^{n}C(F_{U_{r:n}}(x),F_{V_{r:n}}(y)),
\end{eqnarray*}%
where $C(t,s)=\min (t,s)$ and $F_{U_{r:n}}(x)=\dsum\limits_{i=r}^{n}\binom{n%
}{i}x^{i}(1-x)^{n-i}$ and $F_{V_{r:n}}(y)=$ $\dsum\limits_{i=r}^{n}\binom{n}{%
i}y^{i}(1-y)^{n-i}.$

Assume that marginal distributions are $Uniform(0,1)$, i.e. $%
F_{X}(x)=x,0\leq x\leq 1,$ $F_{Y}(y)=y,$ $0\leq y\leq 1.$ \ Denote by $%
(U_{i},V_{i}),i=1,2,...,n$ the random sample from the bivariate distribution 
$C(u,v),$ $0\leq u\leq 1,0\leq v\leq 1$ \ and $U_{1:n}\leq U_{2:n}\leq
\cdots \leq U_{n:n},$ $V_{1:n}\leq V_{2:n}\leq \cdots \leq V_{n:n}$ be
corresponding order statistics. Then from (\ref{x1}) for $r=s$ one has 
\begin{eqnarray*}
P\{U_{r:n} &\leq &u,V_{r:n}\leq v\} \\
&=&\dsum\limits_{i=r}^{n}\sum_{j=r}^{n}%
\sum_{k=a}^{b}c(n,k;i,j)C^{k}(u,v)(u-C(u,v))^{i-k}
\end{eqnarray*}%
\begin{equation}
\times (v-C(u,v))^{j-k}(\bar{C}(u,v))^{n-i-j+k},\text{ }0\leq u\leq 1,0\leq
v\leq 1.  \label{d1}
\end{equation}%
and 
\begin{equation*}
K_{+}^{(n)}(u,v)=\frac{1}{n}\dsum\limits_{r=1}^{n}\dsum\limits_{i=r}^{n}%
\sum_{j=r}^{n}\sum_{k=a}^{b}c(n,k;i,j)C(u,v)^{k}(u-C(u,v))^{i-k}
\end{equation*}%
\begin{equation}
(v-C(u,v))^{j-k}(\bar{C}(u,v))^{n-i-j+k},0\leq u\leq 1,0\leq v\leq 1.
\label{d2}
\end{equation}

\textbf{Lemma 1. }The joint distribution of $U_{r:n}$ and $V_{r:n}$ can be
represented as 
\begin{eqnarray}
P\{U_{r:n} &\leq &u,V_{r:n}\leq v\}  \notag \\
&=&\dsum\limits_{i=r}^{n}\binom{n}{i}C^{i}(u,v)\left[ v-C(u,v)+\bar{C}(u,v)%
\right] ^{n-i}  \notag \\
&&+\dsum\limits_{j=r}^{n}\binom{n}{j}C^{j}(u,v)\left[ u-C(u,v)+\bar{C}(u,v)%
\right] ^{n-j}  \notag \\
&&-\dsum\limits_{i=r}^{n}\binom{n}{i}C^{i}(u,v)(\bar{C}(u,v))^{n-i}  \notag
\\
&&+\dsum\limits_{i=r}^{n}\sum_{j=r}^{n}\sum_{k\neq i\neq
j}c(n,k;i,j)C(u,v)^{k}(u-C(u,v))^{i-k}  \notag \\
&&\times (v-C(u,v))^{j-k}(\bar{C}(u,v))^{n-i-j+k}.  \label{d6}
\end{eqnarray}

\textbf{Proof.} Separating terms for summation in (\ref{d1}) for $%
k=i=j,k=i\neq j$ and $\ k=j\neq i$ we have 
\begin{eqnarray*}
P\{U_{r:n} &\leq &u,V_{r:n}\leq v\} \\
&=&\dsum\limits_{i=r}^{n}c(n,i,i,i)C^{i}(u,v)(\bar{C}(u,v))^{n-i} \\
&&+\dsum\limits_{i=r}^{n}%
\sum_{j=i+1}^{n}c(n,i,i,j)C^{i}(u,v)(u-C(u,v))^{i-i}(v-C(u,v))^{j-i}(\bar{C}%
(u,v))^{n-i-j+i} \\
&&+\dsum\limits_{j=r}^{n}%
\sum_{i=j+1}^{n}c(n,i,i,j)C^{j}(u,v)(u-C(u,v))^{i-j}(v-C(u,v))^{j-j}(\bar{C}%
(u,v))^{n-i-j+j} \\
&&+\dsum\limits_{i=r}^{n}\sum_{j=r}^{n}\sum_{k\neq i\neq
j}c(n,k;i,j)C(u,v)^{k}(u-C(u,v))^{i-k}(v-C(u,v))^{j-k}(\bar{C}%
(u,v))^{n-i-j+k}
\end{eqnarray*}%
\begin{eqnarray}
&=&\dsum\limits_{i=r}^{n}\binom{n}{i}C^{i}(u,v)(\bar{C}(u,v))^{n-i}  \notag
\\
&&+\dsum\limits_{i=r}^{n}\sum_{j=i+1}^{n}\frac{n!}{i!(j-i)!(n-j)!}%
C^{i}(u,v)(v-C(u,v))^{j-i}(\bar{C}(u,v))^{n-j}  \notag \\
&&+\dsum\limits_{j=r}^{n}\sum_{i=j+1}^{n}\frac{n!}{j!(i-j)!(n-i)!}%
C^{j}(u,v)(u-C(u,v))^{i-j}(\bar{C}(u,v))^{n-i}  \notag \\
&&+\dsum\limits_{i=r}^{n}\sum_{j=r}^{n}\sum_{k\neq i\neq
j}c(n,k;i,j)C(u,v)^{k}(u-C(u,v))^{i-k}  \notag \\
&&\times (v-C(u,v))^{j-k}(\bar{C}(u,v))^{n-i-j+k}.  \label{d3}
\end{eqnarray}%
\newline
Consider the second summation in (\ref{d3}) and changing index in the inner
sum as $j-i=k$ \ $(j=i+k,$ $j=i+1\Longrightarrow k=1,j=n\Longrightarrow
k=n-i)$ we have 
\begin{eqnarray}
&&\dsum\limits_{i=r}^{n}\sum_{j=i+1}^{n}\frac{n!}{i!(j-i)!(n-j)!}%
C^{i}(u,v)(v-C(u,v))^{j-i}(\bar{C}(u,v))^{n-j}  \notag \\
&=&\dsum\limits_{i=r}^{n}C^{i}(u,v)\frac{n!}{i!(n-i)!}\sum_{k=1}^{n-i}\frac{%
(n-i)!}{k!(n-i-k)!}(v-C(u,v))^{k}(\bar{C}(u,v))^{n-i-k}  \notag \\
&=&\dsum\limits_{i=r}^{n}\binom{n}{i}C^{i}(u,v)\left[ \sum_{k=0}^{n-i}\frac{%
(n-i)!}{k!(n-i-k)!}(v-C(u,v))^{k}(\bar{C}(u,v))^{n-i-k}-(\bar{C}(u,v))^{n-i}%
\right]  \notag \\
&=&\dsum\limits_{i=r}^{n}\binom{n}{i}C^{i}(u,v)\left[ v-C(u,v)+\bar{C}(u,v)%
\right] ^{n-i}-\dsum\limits_{i=r}^{n}\binom{n}{i}C^{i}(u,v)(\bar{C}%
(u,v))^{n-i}.  \label{d4}
\end{eqnarray}%
Analogously, the third term in (\ref{d3}) can be written as 
\begin{eqnarray}
&&\dsum\limits_{j=r}^{n}\sum_{i=j+1}^{n}\frac{n!}{j!(i-j)!(n-i)!}%
C^{j}(u,v)(u-C(u,v))^{i-j}(\bar{C}(u,v))^{n-i}  \notag \\
&=&\dsum\limits_{j=r}^{n}\binom{n}{j}C^{j}(u,v)\left[ u-C(u,v)+\bar{C}(u,v)%
\right] ^{n-i}  \label{d5} \\
&&-\dsum\limits_{i=r}^{n}\binom{n}{i}C^{i}(u,v)(\bar{C}(u,v))^{n-i}.  \notag
\end{eqnarray}%
Taking into account (\ref{d4}) and (\ref{d5}) in (\ref{d3}) we have (\ref{d6}%
).

\bigskip

\textbf{Theorem 2. }If $C(u,v)=\min (u,v),$ then 
\begin{eqnarray*}
P\{U_{r:n} &\leq &u,V_{r:n}\leq v\} \\
&=&\left\{ 
\begin{array}{ccc}
\dsum\limits_{i=r}^{n}\binom{n}{i}u^{i}\left[ 1-u\right] ^{n-i} & if & u\leq
v \\ 
\dsum\limits_{i=r}^{n}\binom{n}{i}v^{i}\left[ 1-v\right] ^{n-i} & if & u>v%
\end{array}%
\right. \\
&=&\left\{ 
\begin{array}{ccc}
P\{U_{r:n}\leq u\} & if & u\leq v \\ 
P\{V_{r:n}\leq v\} & if & u>v%
\end{array}%
\right.
\end{eqnarray*}

\textbf{Proof. }Consider (\ref{d6}). \ Let $C(u,v)=\min (u,v).$ Then it is
clear that $(u-C(u,v))^{i-k}(v-C(u,v))^{j-k}=0$ for those $k$ satisfying $%
k\neq i$ and $k\neq j,$ because $u-C(u,v)=u-u=0$ if $u\leq v$ and $%
v-C(u,v)=0 $ if $\ u>v.$ Therefore the last term of (\ref{d6}) vanishes and
we have from (\ref{d6}) 
\begin{eqnarray}
P\{U_{r:n} &\leq &u,V_{r:n}\leq v\}  \notag \\
&=&\left\{ 
\begin{array}{ccc}
\dsum\limits_{i=r}^{n}\binom{n}{i}u^{i}\left[ 1-u\right] ^{n-i} & if & u\leq
v \\ 
\dsum\limits_{i=r}^{n}\binom{n}{i}v^{i}\left[ 1-v\right] ^{n-i} & if & u>v%
\end{array}%
\right.  \label{d7} \\
&=&\left\{ 
\begin{array}{ccc}
P\{U_{r:n}\leq u\} & if & u\leq v \\ 
P\{V_{r:n}\leq v\} & if & u>v%
\end{array}%
\right. .  \notag
\end{eqnarray}

\bigskip Consider now, the Baker's Type I BB copula obtained from (\ref{d7}) 
\begin{equation*}
K_{+}^{(n)}(u,v)=\frac{1}{n}\dsum\limits_{r=1}^{n}P\{U_{r:n}\leq
u,V_{r:n}\leq v\}
\end{equation*}%
\begin{equation*}
=\left\{ 
\begin{array}{ccc}
\frac{1}{n}\dsum\limits_{r=1}^{n}\dsum\limits_{i=r}^{n}\binom{n}{i}u^{i}%
\left[ 1-u\right] ^{n-i} & if & u\leq v \\ 
\frac{1}{n}\dsum\limits_{r=1}^{n}\dsum\limits_{i=r}^{n}\binom{n}{i}v^{i}%
\left[ 1-v\right] ^{n-i} & if & u>v%
\end{array}%
\right.
\end{equation*}%
\begin{equation*}
=\left\{ 
\begin{array}{ccc}
u & if & u\leq v \\ 
v & if & u>v%
\end{array}%
\right. =W(u,v)=\min (u,v)
\end{equation*}%
is nothing but the Frech\'{e}t upper bound itself. \ In other words the
Baker's BB distribution with uniform marginals and underlying joint
distribution being Frech\'{e}t upper bound generates the Frech\'{e}t upper
bound.

Now consider the Baker's Type II BB distribution with the "kernel" copula
being Frech\'{e}t upper bound, i.e. 
\begin{equation*}
G_{+}^{(n)}(x,y)=\frac{1}{n}\sum_{r=1}^{n}C(F_{U_{r:n}}(x),F_{V_{r:n}}(y)),
\end{equation*}%
where $C(t,s)=\min (t,s)$ and $F_{U_{r:n}}(x)=\dsum\limits_{i=r}^{n}\binom{n%
}{i}x^{i}(1-x)^{n-i}$ and $F_{V_{r:n}}(y)=$ $\dsum\limits_{i=r}^{n}\binom{n}{%
i}y^{i}(1-y)^{n-i}.$ \ It follows that 
\begin{eqnarray*}
G_{+}^{(n)}(x,y) &=&\frac{1}{n}\sum_{r=1}^{n}\min
(F_{U_{r:n}}(x),F_{V_{r:n}}(y)) \\
&=&\frac{1}{n}\sum_{r=1}^{n}\min (F_{U_{r:n}}(x),F_{V_{r:n}}(y)) \\
&=&\left\{ 
\begin{array}{ccc}
\frac{1}{n}\sum_{r=1}^{n}\dsum\limits_{i=r}^{n}\binom{n}{i}x^{i}(1-x)^{n-i}
& if & x\leq y \\ 
\frac{1}{n}\sum_{r=1}^{n}\dsum\limits_{i=r}^{n}\binom{n}{i}y^{i}(1-y)^{n-i})
& if & x>y%
\end{array}%
\right. \\
&=&\min (x,y).
\end{eqnarray*}%
\bigskip Therefore, if in the constructions $K_{+}^{(n)}(x,y)$ and $%
G_{+}^{(n)}(x,y)$ with $Uniform(0,1)$ marginals one uses the Frech\'{e}t
upper bound copula $W(t,s),$ then obtained copula is again the Frech\'{e}et
upper bound copula. \ From the Example 1 one observes that if in the same
construction one uses FGM copula then different distributions with high
correlation can be obtained. For different constructions satisfying
conditions of the Example 1, we have the following set of inequalities for
correlation coefficient 
\begin{equation*}
\rho _{M}\leq \rho _{G_{-}^{(n)}}\leq \rho _{K_{-}^{(n)}}\leq \rho
_{H_{-}^{(n)}}\leq \rho _{H_{+}^{(n)}}\leq \rho _{K_{+}^{(n)}}\leq \rho
_{G_{+}^{(n)}}\leq \rho _{W},
\end{equation*}%
where \ $\rho _{M}$ is the correlation coefficient of the Frech\'{e}et lower
bound $M(t,s)=\max (t+s-1,0).$

\end{document}